# Causal Inference from Strip-Plot Designs in a Potential Outcomes Framework

Fatemah A. Alquallaf [1*], S. Huda [1] and Rahul Mukerjee [2]

[1] Department of Statistics and OR, Faculty of Science, Kuwait University, P.O. Box-5969, Safat-13060, Kuwait

[2] Indian Institute of Management Calcutta, Joka, Diamond Harbour Road, Kolkata 700104, India

*Corresponding author. E-mail: falquallaf@gmail.com

*Abstract*: Strip-plot designs are very useful when the treatments have a factorial structure and the factors levels are hard-to-change. We develop a randomization-based theory of causal inference from such designs in a potential outcomes framework. For any treatment contrast, an unbiased estimator is proposed, an expression for its sampling variance is worked out, and a conservative estimator of the sampling variance is obtained. This conservative estimator has a nonnegative bias, and becomes unbiased under between-block additivity, a condition milder than Neymannian strict additivity. A minimaxity property of this variance estimator is also established. Simulation results on the coverage of resulting confidence intervals lend support to theoretical considerations.

*Key words*: Between-block additivity; conservative variance estimator; minimaxity; treatment contrast; unbiased estimator.

*2010 Mathematics Subject Classification*: 62K15.

**1. Introduction**

Strip-plot designs are of much practical utility when the treatments have a factorial structure and the factor levels are hard-to-change, so that they have to be applied to larger clusters of experimental units. For instance, in an agricultural field experiment with two factors, irrigation and harvesting, both requiring larger plots, the experimental units can be laid out in several blocks, each block being a rectangular array of rows and columns. One can then randomize the methods of irrigation among the rows and the methods of harvesting among the columns, in each block. A second example where a strip-plot design can be very useful comes from industry. In a semiconductor experiment, it makes much sense to randomize the implant and annealing factors, respectively, among the rows and columns of a rectangular array of wafers, with each such array constituting a block. Again in a painting experiment, primers and paints may be randomized, respectively, among the horizontal and vertical strips on a wall and each such wall may constitute a block. See Casella (2008, Section 5.6.1) and Milliken and Johnson (2009, Chapters 5, 25) for more details and examples.

   With reference to a potential outcomes framework, the present paper aims at developing a randomization-based theory of causal inference from strip-plot designs. In contrast to the traditional method of analysis of such designs, our approach does not hinge on rigid linear model assumptions. After introducing the setup in Section 2, we present the main theoretical results in Section 3 where an unbiased estimator is proposed for any treatment contrast and an expression for its sampling variance is worked out in terms of the various mean squares and products of potential outcomes that arise in a strip-plot setting. We also obtain in Section 3 a conservative estimator of the sampling variance. This conservative estimator is seen to have a nonnegative bias which vanishes under between-block additivity, a condition milder than Neymannian strict additivity (Neymann, 1923/1990). A minimaxity property of our variance estimator, over a class of nonnegative quadratic estimators, is established via a matrix analysis. Finally, Section 4 reports simulation results on the coverage of confidence intervals based on the conservative variance estimator. These are found to match what one expects from theoretical considerations. All longer proofs appear in the appendix.

Before concluding the introduction, we note that causal inference in a potential outcome framework has been of much current interest; recent literature in this general area includes Imbens and Rubin (2015), Dasgupta, Pillai and Rubin (2015) and Mukerjee, Dasgupta and Rubin (2018), where further references are available. Very recently, Zhao, Ding, Mukerjee and Dasgupta (2018) explored at length randomization-based results on causal inference from split-plot designs where the levels of one factor are hard to change while those of the other factor are not. Here, on the other hand, the levels of both factors are hard to change and, consequently, the randomization as well as the final results, including the variance formula, are different from theirs.

**2. Preliminaries**

Consider two experimental factors $F$ and $G$ with levels coded as $1, \ldots, P$ and $1, \ldots, Q$, respectively; $P, Q \geq 2$. The $PQ$ treatment combinations are then $pq$, $1 \leq p \leq P$, $1 \leq q \leq Q$. The factor levels are randomized among $BPQ$ experimental units arranged in $B$ ($\geq 2$) blocks, where each block has $PQ$ units in the form of an $P \times Q$ array. Specifically, we consider a strip-plot design as given by the following randomization which is performed independently for each block and independently for the rows and columns within each block:

(i) levels $1, \ldots, P$ of $F$ are permuted randomly over the $P$ rows, with one row assigned to each level, all such permutations of the $P$ levels over the $P$ rows being equiprobable;

(i) levels $1, \ldots, Q$ of $G$ are permuted randomly over the $Q$ columns, with one column assigned to each level, all such permutations of the $Q$ levels over the $Q$ columns being equiprobable.

Thus, in every block, one unit is assigned to any treatment combination $pq$; this us the unit whose row and column are assigned to levels $p$ and $q$ of $F$ and $G$, respectively.

We next introduce the potential outcomes framework underlying our randomization-based approach. Let $Y_b(rc; pq)$ be the potential outcome from the unit at row $r$ and column $c$ in block $b$, when exposed to treatment combination $pq$, where $1 \leq b \leq B$, $1 \leq r, p \leq P$ and $1 \leq c, q \leq Q$. Then a typical unit-level treatment contrast is of the form $\tau_b(rc) = \Sigma_{p=1}^{P} \Sigma_{q=1}^{Q} l(pq) Y_b(rc; pq)$, where $l(pq)$, $1 \leq p \leq P$, $1 \leq q \leq Q$, are known constants, not all zeros, that sum to zero. Let

$$\bar{Y}_b(pq) = \sum_{r=1}^{P} \sum_{c=1}^{Q} Y_b(rc; pq)/(PQ) \quad \text{and} \quad \bar{Y}(pq) = \sum_{b=1}^{B} \bar{Y}_b(pq)/B \tag{1}$$

denote mean potential outcomes for treatment combination *pq* over the *PQ* units in block *b*, and over all the *BPQ* units, respectively. Then

$$\bar{\tau}_b = \sum_{r=1}^{P} \sum_{c=1}^{Q} \tau_b(rc)/(PQ) = \sum_{p=1}^{P} \sum_{q=1}^{Q} l(pq)\bar{Y}_b(pq) \tag{2}$$

defines a block-level treatment contrast for any fixed *b*, while

$$\bar{\tau} = \sum_{b=1}^{B} \bar{\tau}_b / B = \sum_{p=1}^{P} \sum_{q=1}^{Q} l(pq)\bar{Y}(pq) \tag{3}$$

defines a treatment contrast for the finite population of *BPQ* units.

We shall be concerned with inference on the population level treatment contrast $\bar{\tau}$ on the basis of the potential outcomes observed from the strip-plot design. Note that our definition of $\bar{\tau}$ is quite general and theory developed here will apply immediately to factorial main effect and interaction contrasts as given by specific patterns of the *l(pq)* in $\bar{\tau}$.

## 3. Main results

*3.1 Unbiased estimator of treatment contrast*

For any block *b*, let $Y_b^{\text{obs}}(pq)$ denote the outcome observed from the unit assigned to the treatment combination *pq*. Define

$$\hat{\bar{\tau}}_b = \sum_{p=1}^{P} \sum_{q=1}^{Q} l(pq) Y_b^{\text{obs}}(pq), \ 1 \leq b \leq B, \quad \text{and} \quad \hat{\bar{\tau}} = \sum_{b=1}^{B} \hat{\bar{\tau}}_b / B. \tag{4}$$

Under the strip-plot design, $Y_b^{\text{obs}}(pq)$ can equal any of $Y_b(rc; pq)$, $1 \leq r \leq P$, $1 \leq c \leq Q$, each with probability $1/(PQ)$. Hence by (1),

$$E\{Y_b^{\text{obs}}(pq)\} = \bar{Y}_b(pq), \tag{5}$$

and from (2)-(4), Proposition 1 below is evident. Part (b) of this proposition yields $\hat{\bar{\tau}}$ as an unbiased estimator of the treatment contrast $\bar{\tau}$.

**Proposition 1**. (a) $E(\hat{\bar{\tau}}_b) = \bar{\tau}_b$, $1 \leq b \leq B$, (b) $E(\hat{\bar{\tau}}) = \bar{\tau}$.

*3.2 Sampling variance*

We now proceed to find the sampling variance of the unbiased treatment contrast estimator $\hat{\bar{\tau}}$. As randomization is done independently across blocks, by (4),

$$\text{var}(\hat{\bar{\tau}}) = \sum_{b=1}^{B} \text{var}(\hat{\bar{\tau}}_b)/B^2 = \sum_{b=1}^{B} \sum_{p_1=1}^{P} \sum_{q_1=1}^{Q} \sum_{p_2=1}^{P} \sum_{q_2=1}^{Q} l(p_1 q_1) l(p_2 q_2) \text{cov}\{Y_b^{\text{obs}}(p_1 q_1), Y_b^{\text{obs}}(p_2 q_2)\}/B^2.$$
(6)

Thus, it will suffice to find an expression for $\text{cov}\{Y_b^{\text{obs}}(p_1 q_1), Y_b^{\text{obs}}(p_2 q_2)\}$, for then $\text{var}(\hat{\bar{\tau}})$ can be readily obtained for any treatment contrast $\bar{\tau}$. Consideration of mean squares and products, natural to a strip-plot setting, enables us to obtain such an expression in a neat form in Theorem 1 below.

To that effect, for any fixed $b$, $p_1$, $p_2$, $q_1$ and $q_2$, let

$$M_b^{\text{row}}(p_1q_1, p_2q_2) = Q \sum_{r=1}^{P} \{\bar{Y}_b(r\bullet; p_1q_1) - \bar{Y}_b(p_1q_1)\} \{\bar{Y}_b(r\bullet; p_2q_2) - \bar{Y}_b(p_2q_2)\}/(P-1), \quad (7)$$

$$M_b^{\text{col}}(p_1q_1, p_2q_2) = P \sum_{c=1}^{Q} \{\bar{Y}_b(\bullet c; p_1q_1) - \bar{Y}_b(p_1q_1)\} \{\bar{Y}_b(\bullet c; p_2q_2) - \bar{Y}_b(p_2q_2)\}/(Q-1), \quad (8)$$

$$M_b^{\text{rowcol}}(p_1q_1, p_2q_2) = \sum_{r=1}^{P}\sum_{c=1}^{Q} \{Y_b(rc; p_1q_1) - \bar{Y}_b(r\bullet; p_1q_1) - \bar{Y}_b(\bullet c; p_1q_1) + \bar{Y}_b(p_1q_1)\}$$
$$\times \{Y_b(rc; p_2q_2) - \bar{Y}_b(r\bullet; p_2q_2) - \bar{Y}_b(\bullet c; p_2q_2) + \bar{Y}_b(p_2q_2)\}/\{(P-1)(Q-1)\}, \quad (9)$$

where $\bar{Y}_b(r\bullet; pq) = \Sigma_{c=1}^{Q} Y_b(rc; pq)/Q$, $\bar{Y}_b(\bullet c; pq) = \Sigma_{r=1}^{P} Y_b(rc; pq)/P$, and $\bar{Y}_b(pq)$ is as in (1). We also write $\delta(.,.)$ for Kronecker delta which equals 1 when the two arguments are equal, and 0 otherwise. Then the following result, proved in the appendix, holds.

**Theorem 1**. *For $1 \leq b \leq B$, $1 \leq p_1, p_2 \leq P$ and $1 \leq q_1, q_2 \leq Q$,*

$$\text{cov}\{Y_b^{\text{obs}}(p_1q_1), Y_b^{\text{obs}}(p_2q_2)\}$$
$$= [\{P\delta(p_1,p_2)-1\}M_b^{\text{row}}(p_1q_1, p_2q_2) + \{Q\delta(q_1,q_2)-1\}M_b^{\text{col}}(p_1q_1, p_2q_2)$$
$$+ \{P\delta(p_1,p_2)-1\}\{Q\delta(q_1,q_2)-1\}M_b^{\text{rowcol}}(p_1q_1, p_2q_2)]/(PQ).$$

The proof of Theorem 1 requires much more effort than that of Proposition 1. As seen in the appendix, a lemma, inspired by finite population sampling, and a conditioning argument help.

*3.3 Conservative estimator of sampling variance*

By (6)-(9) and Theorem 2, $\text{var}(\hat{\bar{\tau}})$ involves product terms like $Y_b(rc; p_1q_1)Y_b(rc; p_2q_2)$. Hence, as in other situations with a potential outcome framework (see e.g., Dasgupta, Pillai and Rubin, 2015), unbiased estimation of $\text{var}(\hat{\bar{\tau}})$ is not possible as each experimental unit receives only one treatment combination and, therefore, both $Y_b(rc; p_1q_1)$ and $Y_b(rc; p_2q_2)$ are not observable whenever $p_1q_1 \neq p_2q_2$. On the other hand, $\hat{\bar{\tau}}$ is the arithmetic mean of $\hat{\bar{\tau}}_1, ..., \hat{\bar{\tau}}_B$ by (4), and this prompts one to consider an estimator of $\text{var}(\hat{\bar{\tau}})$ as given by

$$\hat{\text{var}}_0(\hat{\bar{\tau}}) = \sum_{b=1}^{B} (\hat{\bar{\tau}}_b - \hat{\bar{\tau}})^2/\{B(B-1)\}. \quad (10)$$

**Theorem 2**. $E\{\hat{\text{var}}_0(\hat{\bar{\tau}})\} = \text{var}(\hat{\bar{\tau}}) + \Delta_0$, *where* $\Delta_0 = \Sigma_{b=1}^{B} (\bar{\tau}_b - \bar{\tau})^2/\{B(B-1)\}$.

*Proof.* By (4), $\Sigma_{b=1}^{B}(\hat{\bar{\tau}}_b - \hat{\bar{\tau}})^2 = \Sigma_{b=1}^{B}\hat{\bar{\tau}}_b^2 - B\hat{\bar{\tau}}^2$. Hence the result follows noting that

$$E\{\sum_{b=1}^{B}(\hat{\bar{\tau}}_b - \hat{\bar{\tau}})^2\} = \sum_{b=1}^{B}\{\text{var}(\hat{\bar{\tau}}_b) + \bar{\tau}_b^2\} - B\{\text{var}(\hat{\bar{\tau}}) + \bar{\tau}^2\} = B(B-1)\text{var}(\hat{\bar{\tau}}) + \sum_{b=1}^{B}(\bar{\tau}_b - \bar{\tau})^2,$$

using (3), (6) and Proposition 1. □

In view of Theorem 2, $\hat{var}_0(\hat{\bar{\tau}})$ is a conservative estimator of the sampling variance. It has a nonnegative bias $\Delta_0$ which vanishes when $\bar{\tau}_1 = ... = \bar{\tau}_B$. From (2), it is not hard to see that this happens for every treatment contrast if and only if the potential outcomes satisfy a between-block additivity condition as given by the constancy of $\bar{Y}_b(p_1q_1) - \bar{Y}_b(p_2q_2)$ over $1 \leq b \leq B$, for every pair of treatment combinations $p_1q_1$ and $p_2q_2$. This condition is comparable to that of between-whole-plot additivity obtained by Zhao, Ding, Mukerjee and Dasgupta (2018) for split-plot designs though, because of (4), it arises more naturally in our context. Note that this between-block additivity is weaker than Neymannian strict additivity which demands, for every pair of treatment combinations $p_1q_1$ and $p_2q_2$, the constancy of $Y_b(rc; p_1q_1) - Y_b(rc; p_2q_2)$, over $b$, $r$ and $c$,

*3.4 Minimaxity of the conservative variance estimator*

We now establish, in Theorems 3 and 4 below, a minimaxity property of $\hat{var}_0(\hat{\bar{\tau}})$, with regard to its bias $\Delta_0$, over a class of comparable variance estimators. Theorem 3, which is the main plank of this subsection, gives a characterization for the variance estimators in the competing class. Its proof in the appendix involves a matrix analysis and makes use of the sampling variance formula obtained in (6) and Theorem 1. The characterization derived in Theorem 3 paves the way for the final result, Theorem 4, via standard eigenvalue considerations.

Let $Y^{obs} = ((Y_1^{obs})', ...., (Y_B^{obs})')'$ be the $BPQ$x1 vector of observed potential outcomes, where the primes stand for transpose, and each $Y_b^{obs}$ is $PQ$x1 with elements $Y_b^{obs}(pq)$, $1 \leq p \leq P$ and $1 \leq q \leq Q$. Consider now a class, $V$, of variance estimators of the form

$$\hat{var}(\hat{\bar{\tau}}) = (Y^{obs})' A (Y^{obs}), \qquad (11)$$

where $A$ is any known nonnegative definite (nnd) matrix of order $BPQ$, such that $E\{\hat{var}(\hat{\bar{\tau}})\} = var(\hat{\bar{\tau}})$, whenever between-block additivity holds. Thus $V$ consists of nonnegative variance estimators which are quadratic in the potential outcomes, and which, like $\hat{var}_0(\hat{\bar{\tau}})$, become unbiased under between-block additivity. Indeed, from (4), (10) and (11), one can check that $\hat{var}_0(\hat{\bar{\tau}})$ belongs to $V$; cf. (12) below. Suppose the range of each potential outcome is an interval $T$, where $T$ can, in particular, be the real line or the positive part thereof. Then the following result, giving a characterization for the variance estimators in $V$, holds. Here $\hat{\bar{\tau}}_{vec} = (\hat{\bar{\tau}}_1, ..., \hat{\bar{\tau}}_B)'$ and $\bar{\tau}_{vec} = (\bar{\tau}_1, ..., \bar{\tau}_B)'$.

**Theorem 3**. (a) *Every variance estimator in the class V is of the form* $\hat{var}(\hat{\bar{\tau}}) = \hat{\bar{\tau}}'_{vec} U \hat{\bar{\tau}}_{vec}$, *where U is nnd of order B having each diagonal element* $1/B^2$ *and each row sum zero.*

(b) *With* $\hat{var}(\hat{\bar{\tau}})$ *as in* (a), $E\{\hat{var}(\hat{\bar{\tau}})\} = var(\hat{\bar{\tau}}) + \Delta$, *where* $\Delta = \bar{\tau}'_{vec} U \bar{\tau}_{vec}$.

Because $U$ is nnd, by Theorem 3, every variance estimator in $V$ turns out to be conservative, like our $\text{vâr}_0(\hat{\bar{\tau}})$, in the sense of having a nonnegative bias $\Delta$. Note also that by (10) and Theorem 2, $\text{vâr}_0(\hat{\bar{\tau}})$ and its bias $\Delta_0$ are also as per Theorem 3, with

$$U = \{B(B-1)\}^{-1}(I - B^{-1}J) = U_0, \quad (12)$$

say, where $I$ is the identity matrix and $J$ the matrix of ones, both of order $B$. A good estimator in $V$ should keep the bias $\Delta = \bar{\tau}'_{\text{vec}} U \bar{\tau}_{\text{vec}}$ under control. As $\bar{\tau}_{\text{vec}}$ is unknown, in the spirit of the $E$-optimality criterion in experimental design, it makes sense to look for a minimax variance estimator in $V$ that minimizes the maximum of $\Delta$ over spherical contours $\{\bar{\tau}_{\text{vec}} : \bar{\tau}'_{\text{vec}} \bar{\tau}_{\text{vec}} = \rho^2\}$, for every $\rho > 0$. This amounts to finding $U$ in Theorem 3(a), so as to minimize its maximum eigenvalue, denoted by $\lambda_{\max}(U)$. Because $U$ has trace $1/B$ and one eigenvalue zero, arguing as in Mukerjee, Dasgupta and Rubin (2018), $\lambda_{\max}(U) \geq 1/\{B(B-1)\}$, the equality being attained if and only if $U = U_0$. Hence Theorem 3 leads to the following result on the minimaxity of our variance estimator.

**Theorem 4**. *The variance estimator* $\text{vâr}_0(\hat{\bar{\tau}})$ *is uniquely minimax in the class $V$ in the sense of minimizing the maximum bias over* $\{\bar{\tau}_{\text{vec}} : \bar{\tau}'_{\text{vec}} \bar{\tau}_{\text{vec}} = \rho^2\}$, *for every* $\rho > 0$.

## 4. Simulation results

We now present simulation results on the coverage of the confidence interval $\hat{\bar{\tau}} \mp z\{\text{vâr}_0(\hat{\bar{\tau}})\}^{1/2}$ for $\bar{\tau}$, where $z$ is the $(1-\alpha/2)$th quantile of the standard normal distribution, i.e., the target coverage is $1-\alpha$. In the simulation, factors $F$ and $G$ have two and three levels, respectively, and each block is a 2x3 array of units, i.e., $P = R = 2$, $Q = C = 3$. The potential outcomes $Y_b(rc; pq)$ are generated as

$$Y_b(rc; pq) = b + b^h\{\psi(pq) + \xi_b(rc; pq) - \bar{\xi}_b(pq)\}, \quad (13)$$

where $h \geq 0$, $\psi(pq) = \exp\{\frac{1}{2}(p-1.5) + \frac{1}{3}(q-2) + (p-1.5)(q-2)\}$, $\bar{\xi}_b(pq) = \Sigma_{r=1}^2 \Sigma_{c=1}^3 \xi_b(rc; pq)/6$, and $\xi_b(rc; pq)$ ($1 \leq b \leq B$; $r, p = 1, 2$; $c, q = 1, 2, 3$) are iid, each uniform $[-1, 1]$. Then $\bar{Y}_b(pq) = b + b^h \psi(pq)$, so that $\bar{\tau}_b = b^h \psi$, where $\psi = \Sigma_{p=1}^2 \Sigma_{q=1}^3 l(pq)\psi(pq)$. As a result, by (3) and Theorem 2, $\bar{\tau} = b_{h0}\psi$, where $b_{h0} = \Sigma_{b=1}^B b^h / B$, and our variance estimator $\text{vâr}_0(\hat{\bar{\tau}})$ has bias

$$\Delta_0 = \psi^2 \sum_{b=1}^B (b^h - b_{h0})^2 / \{B(B-1)\}. \quad (14)$$

In the simulations reported here, the target coverage is taken as $1-\alpha = 0.95$, and we work with $h = 0, 0.5$, and $B = 20, 40, 60$. By (14), the bias $\Delta_0$ vanishes for $h = 0$, while it is of order $O(1)$ for $h = 0.5$. Therefore, for $h = 0$, the confidence interval $\hat{\bar{\tau}} \mp z\{\text{vâr}_0(\hat{\bar{\tau}})\}^{1/2}$ for $\bar{\tau}$ is expected to attain the target coverage at least for large $B$, whereas for $h = 0.5$, over-coverage is anticipated even for large $B$ unless $\psi^2$ is small. We consider five normalized treatment contrasts in our simulation, namely, $\bar{\tau}^{(s)} = (l^{(s)})' \bar{Y}$, $1 \leq s \leq 5$, where $\bar{Y} = (\bar{Y}(11), \bar{Y}(12), \bar{Y}(13), \bar{Y}(21), \bar{Y}(22), \bar{Y}(23))'$, and

$$l^{(1)} = \frac{1}{\sqrt{6}}(1, 1, 1, -1, -1, -1)', \quad l^{(2)} = \frac{1}{2}(1, 0, -1, 1, 0, -1)', \quad l^{(3)} = \frac{1}{\sqrt{12}}(1, -2, 1, 1, -2, 1)',$$

$$l^{(4)} = \frac{1}{2}(1, 0, -1, -1, 0, 1)', \quad l^{(5)} = \frac{1}{\sqrt{12}}(1, -2, 1, -1, 2, -1)'.$$

Of these, the first one represents main effect $F$, the next two represent man effect $G$, and the last two represent interaction $F \times G$.

For each $h$ and $B$ as stated above, we generate the potential outcomes following (13), and obtain 10000 strip-plot designs via appropriate randomization. For each such design and each of the treatment contrasts mentioned above, the confidence interval $\hat{\bar{\tau}} \mp z\{\hat{\text{var}}_0(\hat{\bar{\tau}})\}^{1/2}$ is found on the basis of the observed potential outcomes. Thus, for each of these contrasts and each $h$ and $B$, the simulated coverage of this interval is obtained from 10000 randomly obtained strip-plot designs. The results are summarized in Table 1. As expected, the simulated coverage comes close to the target for $h = 0$, i.e., $\Delta_0 = 0$ in (14). On the other hand, if $h = 0.5$ then $\Delta_0 > 0$. For the contrasts $\bar{\tau}^{(1)}, \bar{\tau}^{(2)}$ and $\bar{\tau}^{(4)}$, then over-coverage persists even with large $B$ and, in fact, tends to increase slightly with $B$ because the quantity $\sum_{b=1}^{B}(b^h - b_{h0})^2 / \{B(B-1)\}$ in $\Delta_0$ equals 0.0523, 0.0537, 0.0542 for $B = 20$, 40, 60, respectively, approaching 0.556, as $B \to \infty$. With $h = 0.5$, however, the coverage turns out to be surprisingly good for $\bar{\tau}^{(3)}$ and $\bar{\tau}^{(5)}$, as these latter contrasts correspond to a much smaller $\psi^2$ and hence a much smaller bias $\Delta_0$ than the others. Thus the simulation results are in accordance with what one expects from theoretical considerations. Many additional simulations, not reported here and performed, for example, with other choices of $\psi(pq)$ in (13) and other possible values of $h$, also lead to very similar findings.

Table 1. *Simulated coverage of the confidence interval* $\hat{\bar{\tau}} \mp z\{\hat{\text{var}}_0(\hat{\bar{\tau}})\}^{1/2}$, *with* $1 - \alpha = 0.95$

|  | $\bar{\tau}^{(1)}$ | $\bar{\tau}^{(2)}$ | $\bar{\tau}^{(3)}$ | $\bar{\tau}^{(4)}$ | $\bar{\tau}^{(5)}$ |
|---|---|---|---|---|---|
| $h = 0, B = 20$ | 0.934 | 0.934 | 0.933 | 0.938 | 0.937 |
| $h = 0, B = 40$ | 0.943 | 0.943 | 0.942 | 0.942 | 0.942 |
| $h = 0, B = 60$ | 0.944 | 0.945 | 0.945 | 0.947 | 0.947 |
| $h = 0.5, B = 20$ | 0.964 | 0.967 | 0.940 | 0.978 | 0.938 |
| $h = 0.5, B = 40$ | 0.973 | 0.975 | 0.945 | 0.984 | 0.946 |
| $h = 0.5, B = 60$ | 0.976 | 0.976 | 0.948 | 0.986 | 0.947 |

**Appendix**

The following lemma, in the spirit of finite population sampling, will be used repeatedly in proving Theorem 1 in a unified manner that avoids a somewhat cumbersome separate consideration of the four cases arising from equal or unequal $p_1$ and $p_2$, and equal or unequal $q_1$ and $q_2$.

**Lemma A.1**. *Let $(x_{11}, x_{12}),\ldots, (x_{N1}, x_{N2})$ be N pairs of numbers and let $\sigma = \{\sigma(1),\ldots,\sigma(N)\}$ be a random permutation of $\{1,\ldots, N\}$, all such N! permutations being equiprobable. Then*

$$\text{cov}\{x_{\sigma(k_1)1}, x_{\sigma(k_2)2}\} = \frac{N\delta(k_1,k_2)-1}{N(N-1)} \sum_{i=1}^{N}(x_{i1}-\bar{x}_1)(x_{i2}-\bar{x}_2), \quad 1 \le k_1, k_2 \le N,$$

*where $\bar{x}_j = \Sigma_{i=1}^{N} x_{ij}/N$, $j = 1, 2$.*

**Proof of Theorem 1**. For any fixed block $b$, let $\sigma_1 = \{\sigma_1(1),\ldots,\sigma_1(P)\}$ and $\sigma_2 = \{\sigma_2(1),\ldots,\sigma_2(Q)\}$ be random permutation of $\{1,\ldots, P\}$ and $\{1,\ldots, Q\}$, respectively, such that row $\sigma_1(p)$ is assigned to level $p$ of $F$, $1 \le p \le P$, and column $\sigma_2(q)$ is assigned to level $q$ of $G$, $1 \le q \le Q$. Note that

$$\text{cov}\{Y_b^{\text{obs}}(p_1 q_1), Y_b^{\text{obs}}(p_2 q_2)\} = \text{cov}[E\{Y_b^{\text{obs}}(p_1 q_1) | \sigma_1\}, E\{Y_b^{\text{obs}}(p_2 q_2) | \sigma_1\}]$$
$$+ E[\text{cov}\{Y_b^{\text{obs}}(p_1 q_1), Y_b^{\text{obs}}(p_2 q_2) | \sigma_1\}]. \quad (A.1)$$

Since $E\{Y_b^{\text{obs}}(pq) | \sigma_1\} = \bar{Y}_b(\sigma_1(p)\bullet; pq)$ for any treatment combination $pq$, by (7) and Lemma A.1,

$$\text{cov}[E\{Y_b^{\text{obs}}(p_1 q_1) | \sigma_1\}, E\{Y_b^{\text{obs}}(p_2 q_2) | \sigma_1\}] = \{P\delta(p_1, p_2) - 1\} M_b^{\text{row}}(p_1 q_1, p_2 q_2)/(PQ). \quad (A.2)$$

Next, using Lemma A.1 again,

$$\text{cov}\{Y_b^{\text{obs}}(p_1 q_1), Y_b^{\text{obs}}(p_2 q_2) | \sigma_1\} = \text{cov}\{Y_b(\sigma_1(p_1)\sigma_2(q_1); p_1 q_1), Y_b(\sigma_1(p_2)\sigma_2(q_2); p_2 q_2) | \sigma_1\}$$
$$= \frac{Q\delta(q_1, q_2) - 1}{Q(Q-1)} \sum_{c=1}^{Q} Z_b(c; p_1 q_1) Z_b(c; p_2 q_2)$$

where $Z_b(c; pq) = Y_b(\sigma_1(p)c; pq) - \bar{Y}_b(\sigma_1(p)\bullet; pq)$. Hence invoking Lemma A.1 once more,

$$E[\text{cov}\{Y_b^{\text{obs}}(p_1 q_1), Y_b^{\text{obs}}(p_2 q_2) | \sigma_1\}]$$
$$= \frac{Q\delta(q_1, q_2) - 1}{Q(Q-1)} \sum_{c=1}^{Q} [E\{Z_b(c; p_1 q_1)\} E\{Z_b(c; p_2 q_2)\} + \text{cov}\{Z_b(c; p_1 q_1), Z_b(c; p_2 q_2)\}]$$
$$= \frac{Q\delta(q_1, q_2) - 1}{Q(Q-1)} \sum_{c=1}^{Q} [\{\bar{Y}_b(\bullet c; p_1 q_1) - \bar{Y}_b(p_1 q_1)\}\{\bar{Y}_b(\bullet c; p_2 q_2) - \bar{Y}_b(p_2 q_2)\}$$
$$+ \frac{P\delta(p_1, p_2) - 1}{P(P-1)} \sum_{r=1}^{P} \{Y_b(rc; p_1 q_1) - \bar{Y}_b(r\bullet; p_1 q_1) - \bar{Y}_b(\bullet c; p_1 q_1) + \bar{Y}_b(p_1 q_1)\}$$
$$\times \{Y_b(rc; p_2 q_2) - \bar{Y}_b(r\bullet; p_2 q_2) - \bar{Y}_b(\bullet c; p_2 q_2) + \bar{Y}_b(p_2 q_2)\}].$$

$$= \{Q\delta(q_1, q_2) - 1\}[M_b^{\text{col}}(p_1 q_1, p_2 q_2) + \{P\delta(p_1, p_2) - 1\} M_b^{\text{rowcol}}(p_1 q_1, p_2 q_2)]/(PQ).$$

in view of (8) and (9). The above, when combined with (A.1) and (A.2), yields the result. □

**Proof of Theorem 3**. For ease in presentation, the proof is split into three steps.

Step 1 (expectation of variance estimator): Because randomization is done independently across blocks, by (5) and Theorem 1, $E(Y^{\text{obs}}) = \bar{Y}$ and $\text{cov}(Y^{\text{obs}}) = \text{diag}(W_1,\ldots, W_B)$, where $\bar{Y}$ is $BPQ\times 1$ with elements $\bar{Y}_b(pq)$, and each $W_b$ is square of order $PQ$ with elements $\text{cov}\{Y_b^{\text{obs}}(p_1 q_1), Y_b^{\text{obs}}(p_2 q_2)\}$.

Hence, for any variance estimator vâr($\hat{\bar{\tau}}$) in $V$ of the form (11), if we write the matrix $A$ in partitioned form as $(A_{bb*})$, where $A_{bb*}$, $1 \leq b, b* \leq B$, are square of order $PQ$, then

$$E\{\text{vâr}(\hat{\bar{\tau}})\} = \text{tr}[E\{A(Y^{obs})(Y^{obs})'\}] = \sum_{b=1}^{B} \text{tr}(A_{bb}W_b) + \bar{Y}'A\bar{Y}. \quad (A.3)$$

<u>Step 2</u> (use of between-block additivity): In (i)-(iii) below, we consider three configurations for the potential outcomes, each entailing between-block additivity, and obtain useful identities invoking the fact that $E\{\text{vâr}(\hat{\bar{\tau}})\}$ equals var($\hat{\bar{\tau}}$) under such additivity. In the process, we note that the $PQ$ rows and columns of any $A_{bb*}$ correspond to the treatment combinations and write $a_{bb*}(p_1q_1, p_2q_2)$ for the $(p_1q_1, p_2q_2)$th element of $A_{bb*}$.

(i) Suppose the potential outcomes are all equal to some nonzero constant $t$ in the interval $T$ that represents their common range. Then, $\bar{Y} = t\varepsilon$, where $\varepsilon$ is the $BPQ \times 1$ vector of ones, and by (7)-(9), each covariance in Theorem 1 vanishes, i.e., each $W_b$ equals the null matrix. Hence by (6) and (A.3), var($\hat{\bar{\tau}}$) = 0 and $E\{\text{vâr}(\hat{\bar{\tau}})\} = t^2\varepsilon'A\varepsilon$. In this case, Neymannian strict additivity, and hence between-block additivity hold. As a result, $E\{\text{vâr}(\hat{\bar{\tau}})\} = \text{var}(\hat{\bar{\tau}})$. Because $t (\neq 0)$, this yields

$$\varepsilon'A\varepsilon = 0. \quad (A.4)$$

(ii) Next, suppose $Y_1(11;11) = Y_1(22;11) = t_1$, $Y_1(12;11) = Y_1(21;11) = t_2$, and let all other potential outcomes equal $t_0 = (t_1 + t_2)/2$, where $t_1$ and $t_2$ are distinct numbers in $T$. Then by (7)-(9) and Theorem 1, cov$\{Y_1^{obs}(11), Y_1^{obs}(11)\} = (t_1 - t_2)^2/(PQ)$ (= $\theta$, say) and all other covariances in Theorem 1 vanish. Also, $\bar{Y} = t_0\varepsilon$. Hence, recalling the definition of $W_b$, by (A.3) and (A.4), $E\{\text{vâr}(\hat{\bar{\tau}})\} = \theta a_{11}(11,11)$, while by (6), var($\hat{\bar{\tau}}$) $= l(11)^2\theta/B^2$. Moreover, between-block additivity holds because $\bar{Y} = t_0\varepsilon$. As a result, $E\{\text{vâr}(\hat{\bar{\tau}})\} = \text{var}(\hat{\bar{\tau}})$, i.e., $a_{11}(11,11) = l(11)^2/B^2$, as $\theta \neq 0$. In a similar manner,

$$a_{bb}(pq, pq) = l(pq)^2/B^2, \quad 1 \leq b \leq B, 1 \leq p \leq P, 1 \leq q \leq Q. \quad (A.5)$$

(iii) Finally, for some fixed $p_0q_0$ ($\neq 11$), let $Y_1(11;11) = Y_1(22;11) = Y_1(11; p_0q_0) = Y_1(22; p_0q_0) = t_1$, $Y_1(12;11) = Y_1(21;11) = Y_1(12; p_0q_0) = Y_1(21; p_0q_0) = t_2$, and let all other potential outcomes equal $t_0$, where $t_1$, $t_2$ and $t_0$ are as in (ii) above. Then by (7)-(9) and Theorem 1,

$$\text{cov}\{Y_1^{obs}(11), Y_1^{obs}(11)\} = \text{cov}\{Y_1^{obs}(p_0q_0), Y_1^{obs}(p_0q_0)\} = \theta,$$

$$\text{cov}\{Y_1^{obs}(11), Y_1^{obs}(p_0q_0)\} = \text{cov}\{Y_1^{obs}(p_0q_0), Y_1^{obs}(11)\} = \theta\delta_0,$$

where $\theta$ is as in (ii) above and $\delta_0 = \{P\delta(1, p_0) - 1\}\{Q\delta(1, q_0) - 1\}/\{(P-1)(Q-1)\}$, while all other covariances in Theorem 1 vanish. Also, $\bar{Y} = t_0\varepsilon$. Hence, by (A.3)-(A.5) and the definition of $W_b$,

$$E\{\text{vâr}(\hat{\bar{\tau}})\} = \theta\{a_{11}(11,11) + a_{11}(p_0q_0, p_0q_0) + 2\delta_0 a_{11}(11, p_0q_0)\}$$

$$= (\theta/B^2)\{l(11)^2 + l(p_0q_0)^2\} + 2\theta\delta_0 a_{11}(11, p_0q_0)\}$$

and by (6), $\text{var}(\hat{\bar{\tau}}) = (\theta/B^2)\{l(11)^2 + l(p_0q_0)^2 + 2\delta_0 l(11)l(p_0q_0)\}$. Again, between-block additivity holds as $\bar{Y} = t_0\varepsilon$. So, $E\{\text{vâr}(\hat{\bar{\tau}})\} = \text{var}(\hat{\bar{\tau}})$, i.e., $a_{11}(11, p_0q_0) = l(11)l(p_0q_0)/B^2$, as $\theta\delta_0 \neq 0$. In a similar manner, $a_{bb}(p_1q_1, p_2q_2) = l(p_1q_1)l(p_2q_2)/B^2$, for every $b$, $p_1q_1$ and $p_2q_2$ such that $p_1q_1 \neq p_2q_2$. Combining this with (A.5), we get in matrix notation

$$A_{bb} = (1/B^2)ll', \quad 1 \leq b \leq B, \tag{A.6}$$

where $l$ is $PQ \times 1$ with elements $l(pq)$, $1 \leq p \leq P$, $1 \leq q \leq Q$.

Step 3 (matrix analysis): Because $A$ is nnd, we get $A = SS'$, for some matrix $S$. In conformity with the partitioned form $(A_{bb^*})$ of $A$, partition $S$ as $S = [S_1' \ldots S_B']'$, where each $S_b$, $1 \leq b \leq B$, has $PQ$ rows. Then $A_{bb^*} = S_b S_{b^*}'$, $1 \leq b, b^* \leq B$, so that, in view of (A.6), for each $b$, the columns of $S_b$ are scalar multiples of $l$, i.e., $S_b = lu_b'$ for some vector $u_b$. Thus,

$$A_{bb^*} = u_{bb^*}ll', \quad 1 \leq b, b^* \leq B, \tag{A.7}$$

where $u_{bb^*} = u_b'u_{b^*}$. Because $l$ is nonnull, it also follows from (A.7) that

$$u_{bb^*} = (l'A_{bb^*}l)/(l'l)^2, \quad 1 \leq b, b^* \leq B. \tag{A.8}$$

In particular, by (A.6) and (A.8),

$$u_{bb} = 1/B^2, \quad 1 \leq b \leq B. \tag{A.9}$$

Recalling that $Y^{\text{obs}} = ((Y_1^{\text{obs}})', \ldots, (Y_B^{\text{obs}})')'$, from (4), (11) and (A.7), we now obtain

$$\text{vâr}(\hat{\bar{\tau}}) = \sum_{b=1}^{B}\sum_{b^*=1}^{B}(Y_b^{\text{obs}})'A_{bb^*}(Y_{b^*}^{\text{obs}}) = \sum_{b=1}^{B}\sum_{b^*=1}^{B}u_{bb^*}(Y_b^{\text{obs}})'ll'(Y_{b^*}^{\text{obs}})$$

$$= \sum_{b=1}^{B}\sum_{b^*=1}^{B}u_{bb^*}\hat{\bar{\tau}}_b\hat{\bar{\tau}}_{b^*} = \hat{\bar{\tau}}_{\text{vec}}'U\hat{\bar{\tau}}_{\text{vec}}, \tag{A.10}$$

where $U = (u_{bb^*})$ is square of order $B$. Because $\hat{\bar{\tau}}_1, \ldots, \hat{\bar{\tau}}_B$ are independent due to independent randomization across the blocks, by (6), (A.9), (A.10) and Proposition 1(a),

$$E\{\text{vâr}(\hat{\bar{\tau}})\} = \sum_{b=1}^{B}\sum_{b^*=1}^{B}u_{bb^*}\{\text{cov}(\hat{\bar{\tau}}_b, \hat{\bar{\tau}}_{b^*}) + E(\hat{\bar{\tau}}_b)E(\hat{\bar{\tau}}_{b^*})\} = (1/B^2)\sum_{b=1}^{B}\text{var}(\hat{\bar{\tau}}_b) + \sum_{b=1}^{B}\sum_{b^*=1}^{B}u_{bb^*}\bar{\tau}_b\bar{\tau}_{b^*}$$

$$= \text{var}(\hat{\bar{\tau}}) + \bar{\tau}_{\text{vec}}'U\bar{\tau}_{\text{vec}}. \tag{A.11}$$

Note also that by (A.8), $U = \{1/(l'l)^2\}(I \otimes l')A(I \otimes l)$, where $\otimes$ denotes Kronecker product and, as before, $I$ is the identity matrix of order $B$. Hence $U$ is nnd because so is $A$. Therefore, to complete the proof, in view of (A.9)-(A.11), it remains to show that $U$ has each row sum zero. To that effect, without loss of generality, suppose $l(11) \neq 0$ and consider a configuration of the potential outcomes such that $Y_b(rc; 11) = t_1$ for every $b, r, c$, while $Y_b(rc; pq) = t_2$ for every $b, r, c$ and $pq\ (\neq 11)$, as before $t_1$ and $t_2$ being distinct numbers in $T$. Then $\bar{Y}_b(11) = t_1$ for every $b$, while $\bar{Y}_b(pq) = t_2$ for every $b$ and $pq\ (\neq 11)$, so that between-block additivity holds. As a result, $E\{\hat{\text{var}}(\hat{\bar{\tau}})\} = \text{var}(\hat{\bar{\tau}})$ and by (A.11), $\bar{\tau}'_{\text{vec}} U \bar{\tau}_{\text{vec}} = 0$. Hence it is immediate that the nnd matrix $U$ has each sum zero because the $l(pq)$, $1 \leq p \leq P$, $1 \leq q \leq Q$, sum to zero and, therefore, in this case, each element of $\bar{\tau}_{\text{vec}}$, i.e., each $\bar{\tau}_b$ equals the nonzero constant $l(11)(t_1 - t_2)$. □

**Acknowledgement**: This work was supported and funded by Kuwait University Research Grant No. SS04/17.